\documentclass[preprint,12pt]{elsarticle}

\usepackage{amssymb}
\usepackage{amsmath}
\usepackage{amsthm}
\usepackage{mathptmx}
\newtheorem{definition}{Definition}[section]

\newtheorem{Lemma}{Lemma}[section]
\newtheorem{Theorem}{Theorem}[section]
\newtheorem*{Proof}{Proof}
\newtheorem*{Remark}{Remark}
\journal{Nuclear Physics B}

\begin{document}
	
	\begin{frontmatter}
		
		
		
		\title{Hyers-Ulam-Rassias stability of functional equations with parameters} %% Article title
		
		
		\author{Jing Zhang,	Qi Liu,Yongmo Hu,Linlin Fu*,Yuxin Wang,\\
			Jinyu Xia,John Michael Rassias,Choonkil Park,Yongjin Li} %% Author name
		
		
		\affiliation{organization={School of Mathematics and Physics,},%Department and Organization
			addressline={ Anqing Normal University}, 
			city={An Qing},
			postcode={246133}, 
			country={P.R.China}}
		\affiliation{organization={School of Mathmatics and Physics},%Department and Organization
			addressline={Shaoyuan  University}, 
			city={ Shaoyuan},
			postcode={512005}, 
			country={P.R.China}}
		
		\affiliation{organization={Department of Mathematics and Informatics},%Department and Organization
			addressline={National and Kapodistrian University of Athens}, 
			city={ Attikis},
			state={15342},
			country={Greece}}
		\affiliation{organization={Research Institute for Convergence of Basic Science},%Department and Organization
			addressline={Hanyang  University}, 
			city={ Seoul},
			state={04763},
			country={Korea}}
		\affiliation{organization={Department of Mathematics},%Department and Organization
			addressline={Sun Yat-Sen University}, 
			city={ Guangzhou},
			postcode={510275}, 
			country={P. R. China}}
		%% Abstract
		\begin{abstract}
			This paper explores the Hyers-Ulam stability of generalized Jensen additive and quadratic functional equations in \(\beta\)-homogeneous \(F\)-space, showing that approximately satisfying mappings have a unique exact approximating counterpart within a specific bound. 
		\end{abstract}
		\begin{keyword}
			Hyers-Ulam stability,Generalized Jensen additive mapping,Quadratic funtional equations,\(\beta\)-homogeneous \(F\)-space 
			
		\end{keyword}
		
	\end{frontmatter}
	
		\section{Introduction}
		 In 1940, Ulam [39] raised the question of functional equation stability: under what conditions does an approximate homomorphism have a nearby homomorphism? Hyers [40]later solved this for approximate additive mappings in Banach spaces, with subsequent extensions allowing unbounded Cauchy differences. Stabilities for quadratic equations, including those with involution, were also established. 
		 
	Suppose \(X\) is a real inner product space, and \(f: X \rightarrow \mathbb{R}\) is a solution to the orthogonal Cauchy functional equation \(f(x+y) = f(x) + f(y)\) where \(\langle x, y \rangle = 0\). In light of the Pythagorean theorem, \(f(x) = \|x\|^2\) constitutes a solution to this conditional equation. Naturally, this function does not satisfy the additivity equation in all cases. Hence, the orthogonal Cauchy equation is not equivalent to the classical Cauchy equation across the entire inner product space. This phenomenon might underscore the importance of studying the orthogonal Cauchy equation.
		
		Over the past few decades, the stability of functional equations has been extensively explored by numerous mathematicians (see [19]). D.H. Hyers was the first to examine the stability of the Cauchy equation [13], demonstrating that if \(f\) is a mapping from a normed space \(X\) to a Banach space \(Y\) satisfying \(\|f(x+y) - f(x) - f(y)\| \leq \epsilon\) for some \(\epsilon > 0\), then there exists a unique additive mapping g:  \(X\rightarrow Y\)  such that \(\|f(x) - g(x)\| \leq \epsilon\). 
		
		In [38], Liu et al.studied the stability of orthogonal additivity in \(\beta\)-homogeneous \(F\)-spaces and quasi-Banach spaces , showing that if f is a mapping from an Abelian group X to a \(\beta_{2}\)-homogeneous F-space Y satisfying \[\|f(x+y) - f(x) - f(y)\| \leq \varepsilon \] for all \(x, y \in X \)with \(x \perp y\). 
		
		Fu et al.studied the stability of orthogonally Jensen additive and quadratic functional equations in \(F\)-spaces and quasi-Banach spaces [45], showing that if \(f\) is a mapping from an Abelian group \(X\) to a \(\beta\)-homogeneous \(F\)-space \(Y\) satisfying \[\left\|2f\left(\frac{x+y}{2}\right) - f(x) - f(y)\right\| \leq \varepsilon\]for all \(x, y \in X\) with \(x \perp y\) , they also showed that if \(f\) satisfies \[\left\|2f\left(\frac{x+y}{2}\right) + 2f\left(\frac{x-y}{2}\right) - f(x) - f(y)\right\| \leq \varepsilon\] for all \(x, y \in X\) with \(x \perp y\).
		
		Building on such foundational work, subsequent scholars have conducted in-depth research on the stability of various functional equations under diverse spaces and conditions (see [31-37]).With regard to the stability of \(\beta\)-homogeneous, reference may be made to [41-44].
		
		In recent years, the Ulam stability of diverse mathematical entities, such as functional equations and inequalities, differential, difference and integral equations, has drawn extensive attention from scholars. Researchers have delved into how these constructs maintain stability under perturbations, exploring the resilience of solutions and the existence of near - exact counterparts. This growing interest reflects the significance of Ulam stability in understanding the behavior of mathematical models across various fields, from basic functional relationships to complex dynamic systems described by differential and integral equations.
		
		Notably, two key equations have garnered attention in orthogonal contexts: the orthogonally Jensen additive functional equation \[2 f\left(\frac{x+y}{2}\right)=f(x)+f(y),\quad  x \perp y \]and the orthogonally Jensen quadratic functional equation \[2 f\left(\frac{x+y}{2}\right)+2 f\left(\frac{x-y}{2}\right)=f(x)+f(y),\quad  x \perp y\] which characterize orthogonally Jensen additive and quadratic mappings, respectively.

		R. Ger and J. Sikorska [12] conducted research on the orthogonal stability of the Cauchy functional equation \(f(x+y) = f(x) + f(y)\). Specifically, they revealed that if \(f\) is a function mapping from an orthogonality space \(X\) to a real Banach space \(Y\), and for all \(x, y \in X\) with \(x \perp y\), the inequality \(\|f(x+y) - f(x) - f(y)\| \leq \epsilon\) holds for some \(\epsilon > 0\), then there exists precisely one orthogonally additive mapping \(g: X \rightarrow Y\) such that \(\|f(x) - g(x)\| \leq \frac{16}{3}\epsilon\) for all \(x \in X\).
		
		F. Skof [27] was the first to study the stability of the quadratic equation. He proved that if \(f\) is a mapping from a normed space \(X\) to a Banach space \(Y\) satisfying \(\|f(x+y) + f(x-y) - 2f(x) - 2f(y)\| \leq \epsilon\) for some \(\epsilon > 0\), then there exists a unique quadratic function \(g: X \rightarrow Y\) such that \(\|f(x) - g(x)\| \leq \frac{\epsilon}{2}\). P. W. Cholewa [3] extended Skof's theorem by substituting \(X\) with an abelian group \(G\). Subsequently, S. Czerwik [4] generalized Skof's result in the vein of Hyers-Ulam. A number of mathematicians have carried out extensive investigations into the stability problem of functional equations (see [5, 6, 14, 21-24, 29]).The author presented and examined the Ulam stability issue pertaining to the relative Euler-Lagrange functional equation 
		$$
		f\left(\lambda x+\mu y\right)+f\left(\mu x-\lambda y\right)=\left(\lambda^2+\mu^2\right)\left[f\left(x\right)+f\left(y\right)\right]
		$$
		in the works [11-13].
		
		The orthogonally quadratic equation
		\[
		f(x+y) + f(x-y) = 2f(x) + 2f(y), \quad x \perp y
		\]
		was initially studied by F. Vajzović [30] in the case where \(X\) is a Hilbert space, \(Y\) is the scalar field, \(f\) is continuous, and \(\perp\) denotes the Hilbert space orthogonality. Later, H. Drljević [7], M. Fochi [10], M. Moslehian [17, 18], and G. Szabó [28] generalized this result.
		
		The concept of orthogonality has a long history, and over the past few decades, various extensions of it have been put forward. In particular, the proposal of the notion of orthogonality in normed linear spaces has been the focus of extensive efforts by many mathematicians.
		
		Let us recall the standard definition of orthogonality by Rätz [25]:
		
		\begin{definition}
			Let \( X \) be a real linear space with \( \dim X \geq 2 \), and let \( \perp \) denote a binary relation on \( X \) satisfying the following conditions:
			
			(1) For all \( x \in X \), \( x \perp 0 \) and \( 0 \perp x \);
			
			(2) For any non-zero elements \( x, y \in X \setminus \{0\} \), if \( x \perp y \), then \( x \) and \( y \) are linearly independent;
			
			(3) For any \( x, y \in X \) with \( x \perp y \), it holds that \( \alpha x \perp \beta y \) for all scalars \( \alpha, \beta \in \mathbb{R} \);
			
			(4) For every two-dimensional subspace \( P \subseteq X \), every \( x \in P \), and every \( \lambda \in [0, \infty) \), there exists some \( y \in P \) such that \( x \perp y \) and \( x + y \perp \lambda x - y \).
			
			Such an ordered pair \( (X, \perp) \) is referred to as an orthogonality space.
		\end{definition}
		
		In 2010, Fechner and Sikorska [9] investigated the stability of orthogonality and introduced an alternative definition of orthogonality, stated as follows:
		\begin{definition}
			Let \( X \) be an Abelian group, and let \( \perp \) be a binary relation defined on \( X \) with the following properties:
			
			(1) For any \( x, y \in X \) with \( x \perp y \), the relations \( x \perp -y \), \( -x \perp y \), and \( 2x \perp 2y \) hold;
			
			(2) For each \( x \in X \), there exists some \( y \in X \) such that \( x \perp y \) and \( x + y \perp x - y \).
		\end{definition}
		In the context of functional analysis, it is essential to highlight that all orthogonal spaces, along with normed linear spaces exhibiting isosceles orthogonality, satisfy these conditions. Nevertheless, spaces characterized solely by Pythagorean orthogonality no longer meet these requirements.
		
		Despite the successful execution of numerous studies on stability, the non-linear structure inherent in infinite-dimensional $F$-spaces has constrained the availability of corresponding stability results. The non-linear nature of $F$-spaces holds substantial importance in functional analysis and other mathematical domains. A prime illustration is the $L^p([0,1])$ space with $0 < p < 1$, which constitutes an $F$-space (but not a Banach space) when furnished with the metric $d(f, g) = \int|f(x) - g(x)|^p dx$. For a comprehensive understanding of $F$-spaces, interested readers are directed to the literature $[1, 15]$.
		\begin{definition}
			Let $X$ be a linear space. A non-negative function $\|\cdot\|$ is termed an $F$-norm if it adheres to the subsequent conditions:  
			
			(1) $\|x\| = 0$ precisely when $x = 0$;  
			
			(2)  $\|\lambda x\| = \|x\|$ for all scalars $\lambda$ with $|\lambda| = 1$;  
			
			(3) $\|x + y\| \leq \|x\| + \|y\|$ for all $x, y \in X$;  
			
			(4)  $\left\|\lambda_n x\right\| \to 0$ as $\lambda_n \to 0$;  
			
			(5)  $\left\|\lambda x_n\right\| \to 0$ as $x_n \to 0$;  
			
			(6)  $\left\|\lambda_n x_n\right\| \to 0$ when both $\lambda_n \to 0$ and $x_n \to 0$.  
		\end{definition}
		The pair $(X, \|\cdot\|)$ is then denominated an $F^*$-space. An $F$-space is an $F^*$-space that is complete.  
		
		An $F$-norm is said to be $\beta$-homogeneous (where $\beta > 0$) if $\|t x\| = |t|^\beta \|x\|$ for all $x \in X$ and $t \in \mathcal{C}$ (as per $[11, 26]$). A quasi-norm that is $p$-subadditive is referred to as a $p$-norm (with $0 < p < 1$), signifying it satisfies:  
		$$
		\|x + y\|^p \leq \|x\|^p + \|y\|^p \
		$$  
		for all \( x, y \in X.\)
		
		It is important to note that a $p$-subadditive quasi-norm $\|\cdot\|$ gives rise to an $F$-norm. For detailed background information, the reader is advised to consult [1] and [16].  
		\begin{definition}
			\ [2] A quasi-norm $\|\cdot\|$ on a vector space $X$ over a field $K(\mathbb{R})$ is a mapping $X \to [0, \infty)$ with the following attributes:  
			
			(1)  $\|x\| = 0$ if and only if $x = 0$;  
			
			(2)  $\|a x\| = |a| \|x\|$ for all $a \in \mathbb{R}$ and $x \in X$;
			
			(3)  $\|x + y\| \leq C(\|x\| + \|y\|)$ for all $x, y \in X$,  
			
			where $C \geq 1$ is a constant that is independent of $x$ and $y$. The smallest such $C$ for which condition (3) is true is known as the quasi-norm constant of $(X, \|\cdot\|)$.  
		\end{definition}
		A significant result from nonlocally convex theory is the well - known Aoki - Rolewicz theorem [26], which states that for some $0 < p \leq 1$, every quasi-norm has an equivalent $p$-norm.  
		
		Additional results regarding the stability of functional equations in quasi-Banach spaces can be found in [8, 20]. However, these results become even more fascinating and significant when orthogonality is taken into consideration.

		The purpose of this study is to investigate the Hyers-Ulam-Rassias stability of the following parameterized functional equations:
		$$
		f(x+y)=M[f(x)+f(y)],\quad x \perp y
		$$
		where \( M \) is a real number satisfying \( \left|\frac{M-4}{8 M}\right|^\beta+\frac{1}{8^\beta}<1 \) and \( M \neq \frac{1}{2} \);
		and 
		$$
		f(\lambda x+\lambda y)+f(\lambda x-\lambda y)=2 \lambda^2[f(x)+f(y)],\quad x \perp y
		$$
		where \( \lambda \) is an integer and \( \left(\frac{1}{2 \lambda^2}\right)^\beta+\left(\frac{1}{4 \lambda^2}\right)^\beta \leq 1 \),with \(X\) being an Abelian group and \(Y\) being \(\beta\)-homogeneous \(F\)-space as the underlying space setting.

		\section{Stability of the Generalized  Jensen Additive Functional Equations}
		
	In the investigation of the stability of parameterized functional equations, we present the following lemma, which plays a crucial role in establishing the main results of this paper.
		
		\begin{Lemma}
			Let \( X \) be an Abelian group, and \( Y \) be a \( \beta \)-homogeneous \( F \)-space. For \( \varepsilon \geq 0 \), assume \( f: X \to Y \) be a mapping such that for all \( x, y \in X \) , a constant \( C > 0 \) and \( \omega \) is an arbitrary integer, \(\alpha\) and \(\gamma\) are real numbers, one has
			\[
			\left\lVert f( x) - \alpha f(\omega x) +\gamma f(-\omega x) \right\rVert \leq C .\tag{2.1}
			\]
			
			Let
			\[
			h(x, n) = \left\| f( x) - A_{n} f\left(\omega^{n }x\right) + B_{n} f\left(-\omega^{n }x\right) \right\|
			\]
			\text{where} the initial conditions are \[ A_{n+1} = \alpha A_n + \gamma B_n,  B_{n+1} = \gamma A_n + \alpha B_n,  A_0 = 1,  B_0 = 0, |\alpha+\gamma|^\beta+|\alpha-\gamma|^\beta<1,\beta>0\]
			
			and
			\[
			g_n(x) = A_{n} f\left(\omega^n x\right) - B_{n} f\left(-\omega^n x\right), \quad n \in \mathbb{N}.
			\]
			
			Then we have that
			
			(1)\( \left| h(x, n + 1) - h(x, n) \right| \leq C \left[ \ |A_{n}| ^\beta + \ |B_{n}|^\beta \right] \) and moreover, we have
			\[
			h(x, n) \leq C \left( \sum_{n = 1}^{\infty} \left[\ |A_{n}|^\beta + \ |B_{n}|^\beta \right] + 1 \right) 
			\]
			for all \( n \in \mathbb{N} \).	
			
			(2)\( \{ g_n(x) \}_{n \in \mathbb{N}} \) is a Cauchy sequence for every \( x \in X \). Hence, the mapping \( g: X \to Y \) can be defined as:
			\[
			g(x) := \lim_{n \to \infty} g_n(x)
			\]
			and then we have
			\[
			\left\| f( x) - g( x) \right\| \leq C \left( \sum_{n = 1}^{\infty} \left[\ |A_{n}|^\beta + \ |B_{n}|^\beta \right] + 1 \right) 
			\]
			for all \( x \in X \).
		\end{Lemma}
		\begin{Proof} 
			
			First,we have
			\[
			A_n = \frac{(\alpha + \gamma)^n + (\alpha - \gamma)^n}{2}, \quad B_n = \frac{(\alpha + \gamma)^n - (\alpha - \gamma)^n}{2}
			\]
			on account of\quad\( A_{n+1} = \alpha A_n + \gamma B_n,  B_{n+1} = \gamma A_n + \alpha B_n,  A_0 = 1,  B_0 = 0 .\)
			
			So we can get
			\[
			\sum_{n=0}^{\infty} \left( |A_n|^\beta + |B_n|^\beta \right) = \frac{1}{2^\beta} \sum_{n=0}^{\infty} \left[| (\alpha + \gamma)^{n}+(\alpha - \gamma)^{n}|^\beta +|(\alpha + \gamma)^{n}-(\alpha - \gamma)^{n}|^\beta \right]
			\]
			and then
			\[
			\begin{split}
				&|A_n|^\beta + |B_n|^\beta \\
				&\leq \left( |\alpha + \gamma|^n + |\alpha - \gamma|^n \right)^\beta + \left| |\alpha + \gamma|^n - |\alpha - \gamma|^n \right|^\beta\\
				&\leq\	max\left\{|\alpha + \gamma|^{n\beta}, |\alpha - \gamma|^{n\beta}\right\}
			\end{split}
			\]
			
			Because \( |\alpha + \gamma|^\beta + |\alpha - \gamma|^\beta < 1 \), implies \(\max\left\{|\alpha + \gamma|^{\beta}, |\alpha - \gamma|^{\beta}\right\}< 1.\)
			
			Let \( r = \max\left\{|\alpha + \gamma|^{\beta}, |\alpha - \gamma|^{\beta}\right\} \), so we have that \( |A_n|^\beta + |B_n|^\beta \leq r^n \), then we have 
			\[
			\sum_{n=0}^{\infty} \left( |A_n|^\beta + |B_n|^\beta \right) \leq \sum_{n=0}^{\infty} r^n,\quad r<1
			\]
			so the series 
			\[
			\sum_{n = 1}^{\infty} \left( |A_n|^\beta + |B_n|^\beta \right) 
			\]
			converges.
			
			Then,with the help of (2.1), through a simple estimate we obtain
			\[
			\begin{split}
				&\left\| f( x) - A_{n+1} f\left(\omega^{n + 1}x\right) + B_{n+1} f\left(-\omega^{n + 1}x\right) \right\| \\
				&\leq \left\| f( x) - A_{n} f\left(\omega^{n }x\right) + B_{n} f\left(-\omega^{n +}x\right) \right\| \\
				&+ |A_{n}|^\beta \left\| f\left(\omega^{n }x\right) - \alpha f\left(\omega^{n + 1}x\right) + \gamma f\left(-\omega^{n + 1}x\right) \right\| \\
				&+ |B_{n}|^\beta \left\| f\left(-\omega^{n }x\right) - \alpha f\left(-\omega^{n + 1}x\right) + \gamma f\left(\omega^{n + 1}x\right) \right\| \\
				&\leq \left\| f( x) - A_{n} f\left(\omega^{n }x\right) + B_{n} f\left(-\omega^{n }x\right) \right\| \\
				&+ C \left[ |A_{n}|^\beta + |B_{n}|^\beta \right]
			\end{split}
			\]
			which implies
			\[
			\begin{split}
				&\left\lVert f( x) - A_{n} f\left(\omega^{n + 1}x\right) + B_{n} f\left(-\omega^{n + 1}x\right) \right\rVert \\
				&- \left\lVert f( x) - A_{n} f\left(\omega^{n }x\right) + B_{n} f\left(-\omega^{n }x\right) \right\rVert \\
				&\leq C \left[ |A_{n}|^\beta +| B_{n}|^\beta \right].
			\end{split}
			\]
			
			Now we let
			\[
			h(x, n) = \left\| f( x) - A_{n} f\left(\omega^{n }x\right) + B_{n} f\left(-\omega^{n }x\right) \right\|
			\]
			so we have that
			\[
			\lvert h(x, n + 1) - h(x, n) \rvert \leq C \left( |A_{n}|^\beta +| B_{n}|^\beta \right),
			\]
			and then
			\[
			\begin{split}
				h(x, n) &= \sum_{i = 2}^n \bigl( h(x, i) - h(x, i - 1) \bigr) + h(x, 1) \\
				&\leq C \left[ \sum_{n = 1}^\infty \left( |A_{n}|^\beta + |B_{n}|^\beta \right) + 1 \right].
			\end{split}
			\tag{2.2}
			\]
			
			This means that
			\[
			\begin{split}
				&\left\lVert f( x) - A_{n} f\left(\omega^{n }x\right) + B_{n} f\left(-\omega^{n }x\right) \right\rVert \\
				&\leq C \left[ \sum_{n = 1}^\infty \left[ |A_{n}|^\beta + |B_{n}|^\beta \right] + 1 \right], \quad x \in X.
			\end{split}
			\]
			
			The next step is to prove that for each \( x \in X \) the sequence
			\[
			g_n(x) := A_{n} f\left(\omega^n x\right) - B_{n} f\left(-\omega^n x\right), \quad n \in \mathbb{N}
			\]
			is convergent in \( Y \). Since \( Y \) is complete, it suffices to show that \( (g_n(x))_{n \in \mathbb{N}} \) is a Cauchy sequence for every \( x \in X \). Applying estimate (2.1) twice then we have
			\[
			\begin{split}
				\lVert g_n(x) - g_{n + 1}(x) \rVert &= \left\lVert A_{n} \left( f\left(\omega^n x\right) - \alpha f\left(\omega^{n + 1}x\right) + \gamma f\left(-\omega^{n + 1}x\right) \right) \right. \\
				&\quad\quad\left. - B_{n} \left( f\left(-\omega^n x\right) - \alpha f\left(-\omega^{n + 1}x\right) + \gamma f\left(\omega^{n + 1}x\right) \right) \right\rVert \\
				&\leq C \left[ |A_{n}|^\beta + |B_{n}|^\beta \right]
			\end{split}
			\]
			for each \( n \in \mathbb{N} \), which gives us that \( \{g_n(x)\}_{n \in \mathbb{N}} \) is a Cauchy sequence.
			Hence, the mapping \( g: X \to Y \) can be defined as:
			\[
			g(x) := \lim_{n \to \infty} g_n(x)
			\]
			for all \( x \in X \). Combining with (2.2) we have
			\[
			\left\| f( x) - g( x) \right\| \leq C \left( \sum_{n = 1}^{\infty} \left[ |A_{n}|^\beta + |B_{n}|^\beta \right] + 1 \right)
			\]
			with \( x \in X \). \(\square\)
		\end{Proof}
		\begin{Remark}
		The reason why $\omega$ is restricted to an integer here is that $X$ is an Abelian group. Similarly, in the subsequent theorems,  $\lambda$  is required to be integers for the same reason.
		\end{Remark}
		\begin{definition}
			let \( X \) be an Abelian group.A binary relation \( \perp \subseteq X\times X\) is called a generalized Cauchy orthogonality binary relation  on \( X \),if it satisfies the following properties:
			
			(a)For any \(x \in X\), \(0 \perp x\) and \(x \perp 0\);
			
			(b)There exists an integer  such that if \(x \perp x\), then \( x \perp - x.\)
		\end{definition}

		\begin{Theorem}
			Let \(X\) be an Abelian group, and \(Y\) be a \(\beta\)-homogeneous F-space. For \(\varepsilon \geq 0\), assume \(f: X \rightarrow Y\) be a mapping such that for all \(x, y \in X\),\(M\) is a real number satisfying\(
			\left| \frac{M - 4}{8M} \right|^\beta + \frac{1}{8^\beta} < 1
			\) and \(M\neq \frac{1}{2}\), \ one has \(x \perp y\) implies
			\begin{equation}
				\left\lVert f\left( x + y\right)  - M[f(x)+f(y)] \right\rVert \leq \varepsilon.
				\tag{2.3}
			\end{equation}
			
			Then there exists a mapping \(g: X \rightarrow Y\) such that \(x \perp y\) implies
			\begin{equation}
				g\left( x +  y\right)  = M[g(x)+g(y)]
				\tag{2.4}
			\end{equation}
			
			and we have
			\begin{equation}\tag{2.5}
				\begin{aligned}
					&\left\lVert f(x) - g(x) \right\rVert \\
					&\leq \left( \sum_{n=1}^{\infty} \left[ \left( \frac{2^n + 1}{2 \cdot 4^n} \right)^\beta + \left( \frac{2^n - 1}{2 \cdot 4^n} \right)^\beta \right] + 1 \right) \\
					&\quad\left[ \frac{1}{\left\lvert 8M \right\rvert^\beta} + \frac{1}{\left\lvert -16M^2+8M \right\rvert^\beta} + 1 \right] \varepsilon
				\end{aligned}
			\end{equation}
			for all \( x \in 2 X = \{ 2 x \mid x \in X \} .\)
			
			Moreover, the mapping \( g \) is unique on the set \( 2 X \).
		\end{Theorem}
		
		\begin{Proof}
			For all \( x \in X \), we have \( 0 \perp 0 \), \( 0 \perp x \) and \( x \perp 0 \). So according to \( 0 \perp 0 \) and  (2.3) , we have 
			\[
			\left\lVert f(0)  - 2 Mf(0)  \right\rVert \leq \varepsilon.
			\]
			
			Then we can obtain 
			\begin{equation}
				\|f(0)\| \leq \frac{\varepsilon}{|1 - 2M|^{\beta}} .\tag{2.6}
			\end{equation}
			
			According to \( x\perp x \), we can write
			\[\left\| f( 2x )  - 2Mf(x)  \right\| \leq \varepsilon.\]
			
			According to \( x\perp -x \), 	combined with (2.6), it is obvious that
			\begin{equation}
			\left\lvert| f(0)-Mf( 2x )  - Mf(-2x) 
			| \right\rvert \leq \varepsilon.
			\tag{2.7}
		\end{equation}
		
			It is easy to get
			\begin{equation}
				\left\|f( 2x )  - f(-2x) \right\| \leq\frac{\left| \frac{1}{1 - 2M} \right|^\beta + 1 }{|M|^\beta} \varepsilon.
				\tag{2.8}
			\end{equation}

			So it is easy to get
			\begin{align*}
					&\quad\left\lVert f(x) + \frac{M- 4}{8M}f(2x) + \frac{1}{8}f(-2x) \right\rVert\\
					&=\left\lVert \left[f(x) - \frac{1}{2M}f(2x)\right] + \frac{1}{8}\left[f(2x) + f(-2x)\right] \right\rVert\\
					&\leq \varepsilon + \frac{\left\lvert \frac{1}{1 - 2M} \right\rvert^\beta + 1}{\left\lvert 8M \right\rvert^\beta} \varepsilon\\
					&= \left[ \frac{1}{\left\lvert 8M \right\rvert^\beta} + \frac{1}{\left\lvert -16M^2+8M \right\rvert^\beta} + 1 \right] \varepsilon.
			\end{align*}
			
			According to Lemma 2.1,we can easily get
			
			\begin{align*}
				&\quad\left\lVert f(x) - \frac{2^n + 1}{2 \cdot 4^n} f(2^n x) + \frac{2^n - 1}{2 \cdot 4^n} f(-2^n x) \right\rVert\\
				&\leq \left( \sum_{n=1}^{\infty} \left[ \left( \frac{2^n + 1}{2 \cdot 4^n} \right)^\beta + \left( \frac{2^n - 1}{2 \cdot 4^n} \right)^\beta \right] + 1 \right) \\
				&\quad\left[ \frac{1}{\left\lvert 8M \right\rvert^\beta} + \frac{1}{\left\lvert -16M^2+8M \right\rvert^\beta} + 1 \right] \varepsilon\tag{2.9}
		    \end{align*}
			for \(x \in X\) and \(n \in \mathbb{N}\).
			
			Moreover, for each \( x \in X \) the sequence
			\[g_n(x) :=  \frac{2^n + 1}{2 \cdot 4^n} f(2^n x) - \frac{2^n - 1}{2 \cdot 4^n} f(-2^n x), \quad n \in \mathbb{N}\]
			is convergent in \( Y \).
			
			Hence, the mapping 
			\[g(x) := \lim_{n \to \infty} g_n(x)\]
			for all \( x \in X \). Combining with (2.9) we have 
			\begin{align*}
				&\| f( x) - g( x) \| \\
				&\leq \left( \sum_{n=1}^{\infty} \left[ \left( \frac{2^n + 1}{2 \cdot 4^n} \right)^\beta + \left( \frac{2^n - 1}{2 \cdot 4^n} \right)^\beta \right] + 1 \right) \\
				&\quad\left[ \frac{1}{\left\lvert 8M \right\rvert^\beta} + \frac{1}{\left\lvert -16M^2+8M \right\rvert^\beta} + 1 \right] \varepsilon
			\end{align*}
			 for all \(x \in X\).
			
			In order to prove that \( g \) is orthogonally additive, observe first that for \( x,y \in X \) such that \( x \perp y \) and \( n \in \mathbb{N}, \, n > 1 \), we have:
				\begin{align*}
					&\left\lVert g_n(x+y) - M\left[ g_n(x) + g_n(y) \right] \right\rVert \notag \\
					&= \left\lVert \frac{2^n + 1}{2 \cdot 4^n} f(2^n (x+y)) - \frac{2^n - 1}{2 \cdot 4^n} f(-2^n (x+y)) \right\rVert \notag \\
					&\quad - M\left[ \frac{2^n + 1}{2 \cdot 4^n} f(2^n x) - \frac{2^n - 1}{2 \cdot 4^n} f(-2^n x) \right] \notag \\
					&\quad - M\left[ \frac{2^n + 1}{2 \cdot 4^n} f(2^n y) - \frac{2^n - 1}{2 \cdot 4^n} f(-2^n y) \right] \left\lVert \right. \notag \\
					&= \left\lVert \frac{2^n + 1}{2 \cdot 4^n} \left[ f(2^n (x+y)) - M f(2^n x) - M f(2^n y) \right] \right\rVert \notag \\
					&\quad - \frac{2^n - 1}{2 \cdot 4^n} \left[ f(-2^n (x+y)) - M f(-2^n x) - M f(-2^n y) \right] \left\lVert \right. \notag \\
					&\leq \left( \frac{2^n + 1}{2 \cdot 4^n} \right)^\beta \left\lVert f(2^n (x+y)) - M f(2^n x) - M f(2^n y) \right\rVert \notag \\
					&\quad + \left( \frac{2^n - 1}{2 \cdot 4^n} \right)^\beta \left\lVert f(-2^n (x+y)) - M f(-2^n x) - M f(-2^n y) \right\rVert \notag \\
					&\leq \left[ \left( \frac{2^n + 1}{2 \cdot 4^n} \right)^\beta + \left( \frac{2^n - 1}{2 \cdot 4^n} \right)^\beta \right] \varepsilon.
				\end{align*}
			
			Moreover, letting \(n \to \infty\), we get (2.4).
			
			Now, we show the uniqueness of \( g \). Assuming \( g' \) as another mapping satisfying (2.4)and (2.5) that yields:
			\begin{align*}
				&\| g(x) - g'(x) \|\\ &\leq \| g(x) - f(x) \| + \| g'(x) - f(x) \| \\
				&\leq 2 \left( \sum_{n=1}^{\infty} \left[ \left( \frac{2^n + 1}{2 \cdot 4^n} \right)^\beta + \left( \frac{2^n - 1}{2 \cdot 4^n} \right)^\beta \right] + 1 \right) \\
				&\quad\left[ \frac{1}{\left\lvert 8M \right\rvert^\beta} + \frac{1}{\left\lvert -16M^2+8M \right\rvert^\beta} + 1 \right] \varepsilon
			\end{align*}
			for all \( x \in 2 X \).
			
			On the other hand, the mapping \( g - g' \) satisfies (2.4) and thus, in particular, (2.3) with \(\varepsilon = 0\). By applying (2.9) to \( g - g' \) we see that
			\begin{align*}
			&g(2 x) - g'(2 x) \\
			&=  \left( \frac{2^n + 1}{2 \cdot 4^n} \right) \left[ g\left(2^{n+1} x\right) - g'\left(2^{n+1} x\right) \right] \\
			&- \left( \frac{2^n - 1}{2 \cdot 4^n} \right) \left[ g\left(-2^{n+1} x\right) - g'\left(-2^{n+1} x\right) \right]
			\end{align*}
			and therefore
			\begin{align*}
				&\quad\lVert g(2 x) - g'(2 x) \rVert \\
				&\leq \left( \frac{2^n + 1}{2 \cdot 4^n} \right) ^\beta \lVert g(2^{n + 1} x) - g'(2^{n + 1} x) \rVert + \left( \frac{2^n - 1}{2 \cdot 4^n} \right) ^\beta \lVert g(-2^{n + 1} x) - g'(-2^{n + 1} x) \rVert \\
				&\leq 2\left( \sum_{n=1}^{\infty} \left[ \left( \frac{2^n + 1}{2 \cdot 4^n} \right)^\beta + \left( \frac{2^n - 1}{2 \cdot 4^n} \right)^\beta \right] + 1 \right) \\
				&\quad\left[ \frac{1}{\left\lvert 8M \right\rvert^\beta} + \frac{1}{\left\lvert -16M^2+8M \right\rvert^\beta} + 1 \right] \varepsilon
			\end{align*}
			for all \( x \in X \).
			
			Combing the both inequalities,we can easily get the thesis.\(\square\)
		\end{Proof}

		\section{Stability of the Generalized Quadratic Functional Equations}
		
		 In this section, let \( X \) be an Abelian group and let \( \perp \) be a binary relation defined on \( X \) with the properties:
		\begin{enumerate}
			\item[(a)] for any \( x \in X \), \( 0 \perp x \) and \( x \perp 0 \);
			\item[(b)] if \(x,y\in X\) and  \( x \perp x \), then \( - x \perp - x \), \(  x \perp - x \).
		\end{enumerate}

		\begin{Theorem}
			Let \(X\) be an Abelian group, and \(Y\) be a \(\beta\)-homogeneous F-space. For \(\varepsilon \geq 0\), assume \(f: X \rightarrow Y\) be a mapping such that for all \(x, y \in X\) ,where \(\lambda\)  is an integer and \(
			\left( \frac{1}{2\lambda^2} \right)^\beta + \left( \frac{1}{4\lambda^2} \right)^\beta \leq 1
			\), one has \(x \perp y\) implies
			\begin{equation}
				\left\lVert f\left(\lambda x + \lambda y\right) + f\left(\lambda x - \lambda y\right) - 2 \lambda^2\left[f\left(x\right) + f\left(y\right)\right] \right\rVert \leq \varepsilon.
				\tag{3.1}
			\end{equation}
			
			Then there exists a mapping \(g: X \rightarrow Y\) such that \(x \perp y\) implies
			\begin{equation}
				g\left(\lambda x + \lambda y\right) + g\left(\lambda x - \lambda y\right) = 2 \lambda^2\left[g\left(x\right) + g\left(y\right)\right]
				\tag{3.2}
			\end{equation}
			and
			\begin{equation}\tag{3.3}
				\begin{aligned}
					&\left\lVert f(x) - g(x) \right\rVert \\
					&\leq \left(2 + 2^\beta\right)\left|\frac{1}{8\lambda^2}\right|^\beta\left[\left|\frac{1}{2 - 4\lambda^2}\right|^\beta + 1\right]\\
					&\cdot \left( \sum_{n = 1}^\infty \left[ |\frac{2^n + 1}{4^n \lambda^{2n}}|^\beta + |\frac{2^n - 1}{4^n \lambda^{2n}}|^\beta \right] + 1 \right)\varepsilon
				\end{aligned}
			\end{equation}
			for all \( x \in 2\lambda X = \{ 2\lambda x \mid x \in X \} \).
			
			Moreover, the mapping \( g \) is unique on the set \( 2\lambda X \).
		\end{Theorem}
		
		\begin{Proof}
			For all \( x \in X \), we have \( 0 \perp 0 \), \( 0 \perp x \) and \( x \perp 0 \). So according to \( 0 \perp 0 \) and  (3.1) , we have 
			\[
			\left\lVert f(0) + f(0) - 4\lambda^2 f(0)  \right\rVert \leq \varepsilon.
			\]
			
			Then we can obtain 
			\begin{equation}
				\|f(0)\| \leq \frac{\varepsilon}{|2 - 4\lambda^2|^{\beta}}. \tag{3.4}
			\end{equation}
			
			According to \( x \perp x \), we can write
			\[\left\| f(2\lambda x ) + f( 0) - 4\lambda^2 f(x) \right\| \leq \varepsilon.\]
			
			It is easy to get
			
				\[\left\| f(2\lambda x )- 4\lambda^2 f(x) \right\| \leq\left[\left| \frac{1}{2 - 4\lambda^2} \right|^\beta + 1 \right] \varepsilon.\]

			Since \( x\perp -x \), we have
			\[\left\| f(0 ) + f(2 \lambda x) - 2\lambda^2 [f(x) + f(-x)] \right\| \leq \varepsilon.\]
			combined with (3.4),it is obvious that
			\[\left\|   f(2 \lambda x) - 2\lambda^2 [f(x) + f(-x)] \right\| \leq \left[\left| \frac{1}{2 - 4\lambda^2} \right|^\beta + 1 \right]\varepsilon.\]
			
			According to \( -x \perp -x \), we can write
		\[\left\| f(0 ) + f(-2 \lambda x) - 2\lambda^2 [f(-x) + f(-x)] \right\| \leq \varepsilon.\]
		
		It is easy to get
		\[\left\|   f(-2 \lambda x) - 4\lambda^2  f(-x) \right\| \leq \left[\left| \frac{1}{2 - 4\lambda^2} \right|^\beta + 1 \right]\varepsilon.\]
			
			So it is easy to get
			\begin{align*}\tag{3.5}
				&\lVert f( x) - \frac{3}{8\lambda^2} f(2\lambda x) + \frac{1}{8\lambda^2} f(-2\lambda x) \rVert \\
				&= \left\| -\frac{1}{8\lambda^2} \left[ f(2\lambda x) - 4\lambda^2 f(x) \right] + \frac{1}{8\lambda^2} \left[ f(-2\lambda x) - 4\lambda^2 f(-x) \right] \right. \\
				&\quad \left. - \frac{1}{4\lambda^2} \left[ f(2\lambda x) - 2\lambda^2 f(x) - 2\lambda^2 f(-x) \right] \right\|\\
				&\leq \left|\frac{1}{8\lambda^2}\right|^\beta \left[\left|\frac{1}{2 - 4\lambda^2}\right|^\beta + 1\right]\varepsilon + \left|\frac{1}{8\lambda^2}\right|^\beta \left[\left|\frac{1}{2 - 4\lambda^2}\right|^\beta + 1\right]\varepsilon \\
				&+ \left|\frac{1}{4\lambda^2}\right|^\beta \left[\left|\frac{1}{2 - 4\lambda^2}\right|^\beta + 1\right]\varepsilon\\
				&=  \left(2 + 2^\beta\right)\left|\frac{1}{8\lambda^2}\right|^\beta\left[\left|\frac{1}{2 - 4\lambda^2}\right|^\beta + 1\right]\varepsilon.
			\end{align*}
			
			According to Lemma 2.1,we can easily get
			\[
			\begin{split}
				&\left\lVert  f(x) - \frac{2^n + 1}{4^n \lambda^{2n}} f((2\lambda)^n x) + \frac{2^n - 1}{4^n \lambda^{2n}} f\left(-(2\lambda)^n x\right) \right\rVert \\
				&\leq \left(2 + 2^\beta\right)\left|\frac{1}{8\lambda^2}\right|^\beta\left[\left|\frac{1}{2 - 4\lambda^2}\right|^\beta + 1\right]\\
				&\cdot \left( \sum_{n = 1}^\infty \left[ |\frac{2^n + 1}{4^n \lambda^{2n}}|^\beta + |\frac{2^n - 1}{4^n \lambda^{2n}}|^\beta \right] + 1 \right)\varepsilon
			\end{split}\tag{3.6}
			\]
			for all \(x \in X\) and \(n \in \mathbb{N}\).
			
			Moreover, for each \( x \in X \) the sequence
			\[g_n(x) := \frac{2^n + 1}{4^n \lambda^{2n}} f\left((2\lambda)^{n }x\right) - \frac{2^n -1}{4^n \lambda^{2n}} f\left(-(2\lambda)^{n }x\right), \quad n \in \mathbb{N}\]
			is convergent in \( Y \).
			
			Hence, the mapping 
			\[g(x) := \lim_{n \to \infty} g_n(x)\]
			for all \( x \in X \). Combining with (3.6) we have 
			\begin{align*}
				&\| f( x) - g( x) \| \\
				&\leq \left(2 + 2^\beta\right)\left|\frac{1}{8\lambda^2}\right|^\beta\left[\left|\frac{1}{2 - 4\lambda^2}\right|^\beta + 1\right]\\
				&\cdot \left( \sum_{n = 1}^\infty \left[ |\frac{2^n + 1}{4^n \lambda^{2n}}|^\beta + |\frac{2^n - 1}{4^n \lambda^{2n}}|^\beta \right] + 1 \right)\varepsilon, \quad \forall x \in X.
			\end{align*}
			
			In order to prove that \( g \) is orthogonally additive, observe first that for \( x,y \in X \) such that \( x \perp y \) and \( n \in \mathbb{N}, \, n > 1 \), we have:
		\[
		\begin{split}
			&\quad\left\| g(\lambda x + \lambda y) + g(\lambda x - \lambda y) - 2\lambda^2 \bigl[ g(x) + g(y) \bigr] \right\| \\
			&= \Bigg\| \biggl[ \frac{2^n + 1}{4^n \lambda^{2n}} f\bigl( (2\lambda)^n (\lambda x + \lambda y) \bigr) - \frac{2^n - 1}{4^n \lambda^{2n}} f\bigl( -(2\lambda)^n (\lambda x + \lambda y) \bigr) \biggr] \\
			&\quad + \biggl[ \frac{2^n + 1}{4^n \lambda^{2n}} f\bigl( (2\lambda)^n (\lambda x - \lambda y) \bigr) - \frac{2^n - 1}{4^n \lambda^{2n}} f\bigl( -(2\lambda)^n (\lambda x - \lambda y) \bigr) \biggr] \\
			&\quad - 2\lambda^2 \biggl[ \bigl( \frac{2^n + 1}{4^n \lambda^{2n}} f((2\lambda)^n x) - \frac{2^n - 1}{4^n \lambda^{2n}} f(-(2\lambda)^n x) \bigr) + \bigl( \frac{2^n + 1}{4^n \lambda^{2n}} f((2\lambda)^n y) - \frac{2^n - 1}{4^n \lambda^{2n}} f(-(2\lambda)^n y) \bigr) \biggr] \Bigg\| \\
			&= \Bigg\| \frac{2^n + 1}{4^n \lambda^{2n}} \biggl[ f\bigl( (2\lambda)^n (\lambda x + \lambda y) \bigr) + f\bigl( (2\lambda)^n (\lambda x - \lambda y) \bigr) - 2\lambda^2 \bigl( f((2\lambda)^n x) + f((2\lambda)^n y) \bigr) \biggr] \\
			&\quad - \frac{2^n - 1}{4^n \lambda^{2n}} \biggl[ f\bigl( -(2\lambda)^n (\lambda x + \lambda y) \bigr) + f\bigl( -(2\lambda)^n (\lambda x - \lambda y) \bigr) - 2\lambda^2 \bigl( f(-(2\lambda)^n x) + f(-(2\lambda)^n y) \bigr) \biggr] \Bigg\| \\
			&\leq \biggl| \frac{2^n + 1}{4^n \lambda^{2n}} \biggr|^\beta \left\| f\bigl( (2\lambda)^n (\lambda x + \lambda y) \bigr) + f\bigl( (2\lambda)^n (\lambda x - \lambda y) \bigr) - 2\lambda^2 \bigl( f((2\lambda)^n x) + f((2\lambda)^n y) \bigr) \right\| \\
			&\quad + \biggl| \frac{2^n - 1}{4^n \lambda^{2n}} \biggr|^\beta \left\| f\bigl( -(2\lambda)^n (\lambda x + \lambda y) \bigr) + f\bigl( -(2\lambda)^n (\lambda x - \lambda y) \bigr) - 2\lambda^2 \bigl( f(-(2\lambda)^n x) + f(-(2\lambda)^n y) \bigr) \right\| \\
			&\leq \biggl( \biggl| \frac{2^n + 1}{4^n \lambda^{2n}} \biggr|^\beta + \biggl| \frac{2^n - 1}{4^n \lambda^{2n}} \biggr|^\beta \biggr) \varepsilon 
		\end{split}
		\]
			
			Moreover, let\(n \to \infty\), we get (3.2).
			
			Now, we show the uniqueness of \( g \). Assuming \( g' \) as another mapping satisfying (3.3)and (3.4) that yields:
			\begin{align*}
				&\| g(x) - g'(x) \|\\ &\leq \| g(x) - f(x) \| + \| g'(x) - f(x) \| \\
				&\leq 2 \left(2 + 2^\beta\right)\left|\frac{1}{8\lambda^2}\right|^\beta\left[\left|\frac{1}{2 - 4\lambda^2}\right|^\beta + 1\right]\\
				&\cdot \left( \sum_{n = 1}^\infty \left[ |\frac{2^n + 1}{4^n \lambda^{2n}}|^\beta + |\frac{2^n - 1}{4^n \lambda^{2n}}|^\beta \right] + 1 \right)\varepsilon
			\end{align*}
			for all \( x \in 2\lambda X \).
			
			On the other hand, the mapping \( g - g' \) satisfies (3.2) and thus, in particular, (3.1) with \(\varepsilon = 0\). By applying (3.6) to \( g - g' \) we see that
			\begin{align*}
			&\quad\ g(2\lambda x) - g'(2\lambda x) \\
			&= \frac{2^n + 1}{4^n \lambda^{2n}} \left[ g\left( (2\lambda)^{n + 1} x \right) - g'\left( (2\lambda)^{n + 1} x \right) \right] - \frac{2^n - 1}{4^n \lambda^{2n}} \left[ g\left( -(2\lambda)^{n + 1} x \right) - g'\left( -(2\lambda)^{n + 1} x \right) \right]
		\end{align*}
			and therefore
			\begin{align*}
				&\quad\lVert g(2\lambda x) - g'(2\lambda x) \rVert \\
				&\leq \left| \frac{2^n + 1}{4^n \lambda^{2n}} \right|^\beta \left\| g\left( (2\lambda)^{n + 1} x \right) - g'\left( (2\lambda)^{n + 1} x \right) \right\| + \left| \frac{2^n - 1}{4^n \lambda^{2n}} \right|^\beta \left\| g\left( -(2\lambda)^{n + 1} x \right) - g'\left( -(2\lambda)^{n + 1} x \right) \right\|\\
				&\leq \left(  |\frac{2^n + 1}{4^n \lambda^{2n}}|^\beta +  |\frac{2^n - 1}{4^n \lambda^{2n}}|^\beta \right) \cdot 2 \left(2 + 2^\beta\right)\left|\frac{1}{8\lambda^2}\right|^\beta\left[\left|\frac{1}{2 - 4\lambda^2}\right|^\beta + 1\right]\\
				&\cdot \left( \sum_{n = 1}^\infty \left[ |\frac{2^n + 1}{4^n \lambda^{2n}}|^\beta + |\frac{2^n - 1}{4^n \lambda^{2n}}|^\beta \right] + 1 \right)\varepsilon
			\end{align*}
			for \( x \in X \).
			
			Combing the both inequalities,we can easily get the thesis.\(\square\)

		\end{Proof}
		
		\
		
		\

		\noindent \textbf{Acknowledgments.} Thanks to all the members of the Functional Analysis Research team of the College of Mathematics and Physics of Anqing Normal University for their discussion and correction of the difficulties and errors encountered in this paper.
		
		\noindent \textbf{Funding.} This research work was funded by Anhui Province Higher Education Science Research Project (Natural Science), 2023AH050487.
		
		\noindent \textbf{Conflict of interest.} The authors declare that there is no conflict of interest in publishing the article.
		
		%%       \citet{<label>} ==> Lamport [21]
		%%       \citep{<label>} ==> [21]
		%%
		
		
		

	\end{document}